\documentclass[11pt]{amsart}
\usepackage{a4wide}
\usepackage{xcolor}
\definecolor{darkred}{rgb}{0.4,0,0}
\definecolor{darkgreen}{rgb}{0,0.4,0}
\definecolor{darkblue}{rgb}{0,0,0.4}
\usepackage[colorlinks,citecolor=darkgreen,linkcolor=darkred,urlcolor=darkblue,filecolor=darkred]{hyperref}
\usepackage[noabbrev, nameinlink]{cleveref}
\usepackage{newpxtext}
\AtBeginDocument{%
  \DeclareFontShape{\encodingdefault}{\rmdefault}{m}{sl}{<-> pplro7t}{}%
  \DeclareFontShape{\encodingdefault}{\rmdefault}{b}{sl}{<-> pplbo7t}{}%
}
\usepackage{newpxmath}
\usepackage[style=alphabetic,giveninits=true,isbn=false,maxnames=99,maxalphanames=99]{biblatex}
\DeclareFieldFormat{title}{\mkbibquote{#1}}
\DeclareFieldFormat[misc]{date}{\iffieldequalstr{eprinttype}{arXiv}{\mkbibemph{Preprint} (#1)}{#1}}
\AtEveryBibitem{\ifentrytype{unpublished}{\clearfield{year}}{}}
\DeclareFieldFormat[unpublished]{note}{\mkbibemph{#1}}

\AtEveryBibitem{\clearfield{issn}}
\AtEveryCitekey{\clearfield{isbn}}
\AtEveryBibitem{\clearfield{number}}
\AtEveryBibitem{\ifentrytype{unpublished}{}{\ifentrytype{online}{}{\clearfield{url}}}}
\renewbibmacro{in:}{%
\ifentrytype{article}{}{%
\printtext{\bibstring{in}\intitlepunct}%
}%
}
\usepackage{titlesec}
\titleformat{\section}{\normalfont\large\bfseries\scshape}{\mdseries\thesection.}{1em}{}
\titlespacing*{\section}{0pt}{.7\linespacing plus \linespacing}{.5\linespacing}
\titleformat{\subsection}[runin]{\bfseries}{\mdseries\thesubsection.}{.5em}{}[.]
\titlespacing*{\subsection}{0pt}{.5\linespacing plus .7\linespacing}{.5em}
\titleformat{\subsubsection}[runin]{\itshape}{\upshape\thesubsubsection.}{.5em}{}[.]
\titlespacing*{\subsubsection}{0pt}{.5\linespacing plus .7\linespacing}{.5em}
\makeatletter
\def\@settitle{\begin{center}%
\baselineskip14\p@\relax
\bfseries
\large%<- new
\uppercasenonmath\@title
\@title
\end{center}%
}
\g@addto@macro\bfseries{\boldmath}
\makeatother
\usepackage[textsize=footnotesize]{todonotes}
\usepackage{mathtools} %psmallmatrix
\usepackage[shortlabels]{enumitem}
\newlist{enumcases}{enumerate}{2}
\setlist[enumcases,1]{
label={(\arabic*)},
ref=\arabic*} 
\creflabelformat{enumcasesi}{#2(#1)#3}
\setlist[enumcases,2]{
label={\alph*)},
ref=\theenumcasesi.\alph*} 
\creflabelformat{enumcasesii}{#2(#1)#3} 
\crefname{enumcasesi}{case}{cases}
\Crefname{enumcasesi}{Case}{Cases}
\crefname{enumcasesii}{case}{cases}
\Crefname{enumcasesii}{Case}{Cases}
\makeatletter
\@definecounter{formula}
\let\c@formula\c@equation

\def\endformula{\eqno\hbox{\@eqnnum}$$\@ignoretrue}
\makeatother

\crefname{formula}{formula}{formulas}
\Crefname{formula}{Formula}{Formulas}
\usepackage{caption}
\usepackage[labelfont=rm,labelformat=parens]{subcaption}
\usepackage{tikz}
\usetikzlibrary{decorations.pathmorphing,decorations.pathreplacing,backgrounds,patterns,patterns.meta}

\makeatletter
\tikzset{
dot diameter/.store in=\dot@diameter,
dot diameter=1pt,
dot spacing/.store in=\dot@spacing,
dot spacing=2pt,
dots/.style={
line width=\dot@diameter,
line cap=round,
dash pattern=on 0pt off \dot@spacing
}
}
\makeatother
\def\ww{8}
\def\www{3}
\def\wwww{4}
\def\hh{4}
\def\tt{2}
\def\tzero{90/\ww}
\def\squaresize{1}

\newcommand{\cylinder}{
\filldraw[fill=lightgray!33!white,opacity=.618] (0,{\squaresize*\hh}) arc (-180:180:{\squaresize*\ww/6} and {\squaresize*\ww/22.5});
\foreach\i in {0,...,\hh}{
\draw[thick] (0,{\squaresize*\i}) arc (180:0:{\squaresize*\ww/6} and {\squaresize*\ww/22.5});
}
\foreach\i in {\wwww,...,\ww}{
\draw[thick] (0,0) arc (-180:-\i*360/\ww+\tzero:{\squaresize*\ww/6} and {\squaresize*\ww/22.5}) -- ++(0,{\squaresize*\hh});
}

\fill[opacity=.618] (0,0) arc (180:360:{\squaresize*\ww/6} and {\squaresize*\ww/22.5}) -- ++(0,{\squaresize*\hh}) arc (360:180:{\squaresize*\ww/6} and {\squaresize*\ww/22.5}) -- cycle;
\foreach\i in {1,...,\www}{
\draw[thick] (0,0) arc (-180:-\i*360/\ww+\tzero:{\squaresize*\ww/6} and {\squaresize*\ww/22.5}) -- ++(0,{\squaresize*\hh});
}
\draw[thick] (0,0) arc (-180:-\www*360/\ww+\tzero:{\squaresize*\ww/6} and {\squaresize*\ww/22.5}) node {$\bullet$} coordinate (A) -- ++(0,{\squaresize*\hh});
\draw[thick] (0,0) arc (-180:(-\www+\tt)*360/\ww+\tzero:{\squaresize*\ww/6} and {\squaresize*\ww/22.5}) coordinate (B) -- ++(0,{\squaresize*\hh}) node {$\bullet$};

\draw[very thick] (0,0) arc (180:360:{\squaresize*\ww/6} and {\squaresize*\ww/22.5});
\draw[very thick] (0,0) -- (0,{\squaresize*\hh}) ({2*\squaresize*\ww/6},0) -- ({2*\squaresize*\ww/6},{\squaresize*\hh});
\foreach\i in {0,...,\hh}{
\draw[thick] (0,{\squaresize*\i}) arc (-180:0:{\squaresize*\ww/6} and {\squaresize*\ww/22.5});
}
\draw[very thick] (0,{\squaresize*\hh}) arc (-180:180:{\squaresize*\ww/6} and {\squaresize*\ww/22.5});
}
\newcommand{\cylinderparams}{
\draw[thick] (0,0) arc (-180:(-\www+\tt/2)*360/\ww+\tzero:{\squaresize*\ww/6} and {\squaresize*\ww/22.5}) coordinate (O);
\draw[->, very thick] (O) ++(0,-.5) node[below] {$t$} arc ((-\www+\tt/2)*360/\ww+\tzero:-\www*360/\ww+\tzero:{\squaresize*\ww/6} and {\squaresize*\ww/22.5});
\draw[->,very thick] (O) ++(0,-.5) arc ((-\www+\tt/2)*360/\ww+\tzero:(-\www+\tt)*360/\ww+\tzero:{\squaresize*\ww/6} and {\squaresize*\ww/22.5});
\draw[<->,very thick] ({\squaresize*\ww/6},{\squaresize*\hh-\squaresize*\ww/22.5+.5}) arc (-90:270:{\squaresize*\ww/6} and {\squaresize*\ww/22.5}) node[midway,above] {$w$};
\draw[<->,very thick] (0-.5,0) -- (0-.5,{\squaresize*\hh}) node[midway,left] {$h$};
}

\newcommand{\fundom}{
\begin{scope}
\clip (0,0) -- ({\squaresize*\ww},0) -- ({\squaresize*(\ww+\tt)},{\squaresize*\hh}) -- ({\squaresize*\tt},{\squaresize*\hh}) -- cycle;
\filldraw[opacity=.618] (0,0) -- ({\squaresize*\ww},0) -- ({\squaresize*(\ww+\tt)},{\squaresize*\hh}) -- ({\squaresize*\tt},{\squaresize*\hh}) -- cycle;
\foreach\i in {0,...,\hh}{
\draw[thick] (0,{\squaresize*\i}) -- ({\squaresize*(\ww+\tt)},{\squaresize*\i});
}
\foreach\i in {0,...,\ww}{
\draw[thick] ({\i*\squaresize},0) -- ({\i*\squaresize},{\squaresize*\hh});
}
\foreach\i in {0,...,\tt}{
\draw[thick] ({(\i+\ww)*\squaresize},0) -- ({(\i+\ww)*\squaresize},{\squaresize*\hh});
}
\end{scope}
\draw[very thick] (0,0) -- ({\squaresize*\ww},0) -- ({\squaresize*(\ww+\tt)},{\squaresize*\hh}) -- ({\squaresize*\tt},{\squaresize*\hh}) -- cycle;
}
\newcommand{\fundomparams}{
\draw[<->,very thick] (0,0-.5) -- node[midway,below] {$w$} ({\squaresize*\ww},0-.5);
\draw[dots] (0,0) -- (0,0-.5) ({\squaresize*\ww},0) -- ({\squaresize*\ww},0-.5);
\draw[<->,very thick] (0-.5,0) -- node[midway,left] {$h$} (0-.5,{\squaresize*\hh});
\draw[dots] (0,0) -- (0-.5,0) ({\squaresize*\tt},{\squaresize*\hh}) -- (0-.5,{\squaresize*\hh});
\draw[<->,very thick] (0,{\squaresize*\hh+.5}) -- node[midway,above] {$t$} ({\squaresize*\tt},{\squaresize*\hh+.5});
\draw[dots] (0,0) -- (0,{\squaresize*\hh+.5}) ({\squaresize*\tt},{\squaresize*\hh}) -- ({\squaresize*\tt},{\squaresize*\hh+.5});
}
\newtheorem{theorem}{Theorem}[section]
\newtheorem{proposition}[theorem]{Proposition}
\newtheorem{corollary}[theorem]{Corollary}
\newtheorem{lemma}[theorem]{Lemma}
\theoremstyle{remark}
\newtheorem{remark}[theorem]{Remark}

\providecommand{\defn}[1]{{\fontseries{b}\fontshape{sl}\selectfont #1}}
\newcommand{\NN}{\mathbf{N}}
\newcommand{\ZZ}{\mathbf{Z}}
\newcommand{\RR}{\mathbf{R}}
\newcommand{\CC}{\mathbf{C}}
\newcommand{\TT}{\mathbf{T}}
\newcommand{\params}{\mathcal{P}}

\renewcommand{\H}{\mathcal{H}}
\newcommand{\q}{\mathfrak{q}}
\newcommand{\divides}{\mathrel{|}}
\DeclarePairedDelimiter{\floor}{\lfloor}{\rfloor}
\newcommand{\SL}{\mathrm{SL}}
\newcommand{\slz}{{\SL(2,\ZZ)}}
\newcommand{\diff}{\mathrm{d}}
\begin{document}
% -------------------------------------------------------------------------
\title{Square-tiled tori}
\author{Angel Pardo}
\thanks{This work was partially supported by Centro de Modelamiento Matemático (CMM), ACE210010 and FB210005, BASAL funds for centers of excellence from ANID-Chile, 
and by ANID-Chile through the FONDECYT Postdoctorado 3190257 and FONDECYT Regular 1221934 grants.}
\address{
Departamento de Matemática y Ciencia de la Computación\\
Universidad de Santiago de Chile\\
Las Sophoras 173, Estación Central, Santiago, Chile.}
\email{angel.pardo@usach.cl}
\makeatletter
\@namedef{subjclassname@2020}{%
\textup{2020} Mathematics Subject Classification}
\makeatother
%\date{\today}
%\dedicatory{}
\begin{abstract}
We study square-tiled tori, that is, tori obtained from a finite collection of unit squares by parallel side identifications.
Square-tiled tori can be parametrized in a natural way that allows to count the number of square-tiled tori tiled by a given number of square tiles.
There is a natural $\slz$-action on square-tiled tori and we classify $\slz$-orbits using two numerical invariants that can be easily computed.
We deduce the exact size of every $\slz$-orbit.
In particular, this answers a question by M.~Bolognesi on the number of \emph{cyclic} covers of the torus, which corresponds to particular $\slz$-orbits of square-tiled tori.
We also give the asymptotic behavior of the number of cyclic square-tiled tori.
\end{abstract}
% -------------------------------------------------------------------------
\maketitle
%--------------------------------------------------------------------------
% -------------------------------------------------------------------------
\section{Introduction}

A \defn{square-tiled surface} is an orientable connected surface obtained from a finite collection of unit squares after identifications of pairs of parallel sides by translations in a way that right vertical side of a square can only be glued to a left vertical side of a square and viceversa, and similarly for top and bottom sides of squares.

The fundamental example is the square torus $\TT = \CC/\ZZ[i]$, which is obtained from the unit square from identification by translations of parallel sides, as in \Cref{f:square-torus}.

\begin{figure}[ht]
\begin{tikzpicture}[scale=2]
\filldraw[fill=lightgray,very thick,sloped,midway] (0,0) -- node {$|$} (0,1) -- node {$||$} (1,1)-- node {$|$} (1,0) -- node {$||$} cycle;
\end{tikzpicture}
\caption{The square torus obtained from the unit square by opposite side identifications.}
\label{f:square-torus}
\end{figure}

In this work, we focus on \defn{square-tiled tori}, that is, square-tiled surfaces of genus one.
We expose several tools and ideas from the general theory of square-tiled surfaces, in the simplest case.
One goal of this article is to illustrate and introduce some basics of square-tiled surfaces to the interested reader.
On the other hand, most of the results are already known in other contexts, but have not been properly studied though the lens of square-tiled surfaces.
Another aim is to translate and prove these results in the language of square-tiled surfaces and provide a reference within this context.
Furthermore, all our results are stated and proven in an elementary way.

The starting point of this work is a question posed by M.~Bolognesi (by personal communication to A.~Zorich) regarding the counting of cyclic covers of the torus, after his joint work on cyclic covers of the sphere \cite{Bolognesi-Giansiracusa}.

\subsection*{Translation surfaces}
A square-tiled surface is a particular type of a larger class of surfaces known as \defn{translation surfaces}: orientable connected surfaces that can be obtained by edge-to-edge gluing of polygons in $\RR^2$ using translations only.
There is a one-to-one correspondence between translation surfaces and non-trivial Abelian differentials on Riemann surfaces.
For an introduction and general references to the subject of translation surfaces, we refer the reader to the surveys of Zorich~\cite{Zorich:survey} and Forni--Matheus~\cite{Forni-Matheus:survey}.

It follows that square-tiled surfaces can be identified with pairs $(X, \omega)$, where $X$ is a Riemann surface (complex curve) which is a finite cover of $\TT$, say via $\pi\colon X\to\TT$, branched at most at the origin $0 \in \TT$, and $\omega$ is the Abelian differential $\omega = \pi^*(\diff z)$.

\begin{remark}
There are also other alternative definitions for square-tiled surfaces.
In particular, there is a relevant combinatorial description in terms of permutations.
We refrain from giving this characterization, as it is not of much utility in the genus one case, addressed in this work.
In our case, a more natural combinatorial model is the parametrization given in \Cref{t:parametrization}.%, below.
\end{remark}

Square-tiled surfaces are also of general interest in the study of translation surfaces and their moduli spaces, and in Teichmüller dynamics and geometry.
In fact, they are a source of many interesting examples of diverse phenomena in geometry and dynamical systems (see, e.g., \cite{Eskin-Kontsevich-Zorich:square-tiled,Forni-Matheus-Zorich:square-tiled}).
Moreover, they can be seen as \emph{integer points} in moduli spaces of abelian differentials and, in particular, asymptotics on the number of square-tiled surfaces in each moduli space allows to compute their volume (see, e.g., \cite{Eskin-Okounkov,Zorich:square-tiled}).

\subsection*{Flat geometry}
Any translation surface has a canonical flat metric, the one obtained form $\RR^2$, with conical singularities with angles that are a multiples of $2\pi$ at zeros of the Abelian differential.
Since the holonomy lies in $2\pi\ZZ$, given a tangent vector at a point, there is a straight-line trajectory emanating from that point.
These trajectories correspond to the geodesics for the flat metric.
If such a trajectory is closed and does not pass through a cone point, then there is a set of parallel trajectories homotopic to it and that determine a maximal isometrically immersed flat cylinder with boundary, which we call simply \defn{cylinder}.

In the case of square-tiled surfaces, any rational slope define a \defn{cylinder decomposition} in such direction.
The surface is covered by a finite family of cylinders in that direction with piecewise boundary identifications.

For instance, in the case of the horizontal direction, starting from any unit square tile, one can go horizontally to the square to its right and continue going right until one comes back to the starting square.
This gives an immersed cylinder isometric to $\RR/w\ZZ \times [0,1]$, for some $w \in \NN$ that corresponds to the number of squares one had to go through before coming back to the original one.
One can then extend this cylinder vertically and get a maximal immersed (horizontal) cylinder isometric to $\RR/w\ZZ \times [0,h]$, for some $h \in \NN$ corresponding to the number of times the cylinder can be extended.
The parameters $w,h$ are called, respectively, the \emph{width} and the \emph{height} of such horizontal cylinder.
In order to recover the surface one has to determine all the cylinders in a cylinder decomposition together with the information about how to identify their boundaries together.

\subsection*{Square-tiled tori}
Square-tiled surfaces of genus one have notably simpler descriptions than general square-tiled surfaces.
For example, results connecting the geometry of a surface to its topology ---as the Gauss--Bonnet formula or the Euler--Poincaré formula--- already shows that, in this case, the metric is completely flat, with no conical singularities or zeros of the Abelian differential.
An elementary geometrical way of thinking about this fact is that if one start in any unit square, then go to the square to its right, then up, then left and then down, one returns to the starting square, as in \Cref{f:flat}.
This is never the case in higher genus as necessarily at least one corner would have a cone angle greater than $2 \pi$.

\begin{figure}
\begin{tikzpicture}[scale=2,fill=lightgray,very thick]
\fill[decorate,decoration={zigzag,amplitude=.25mm,segment length=5mm}] (.8,.8) rectangle (-.8,-.8);
\draw[dots] (-.85,-.85) grid (.85,.85);
\draw (-.6,-.6) grid (.6,.6);
\draw[->] (1/6,0.28867513459) arc (60:120:1/3);
\draw[->] (-0.28867513459,1/6) arc (150:210:1/3);
\draw[->] (-1/6,-0.28867513459) arc (240:300:1/3);
\draw[->] (0.28867513459,-1/6) arc (330:390:1/3);
\end{tikzpicture}
\caption{Local picture of a square-tiled tori around every corner.}
\label{f:flat}
\end{figure}

On the other hand, in the genus one case, any (maximal) cylinder necessarily covers the entire surface.
That is, any cylinder decomposition consists in only one cylinder with its boundary identified by translation.
In particular, if one considers the width $w$ and the height $h$ of \emph{the} horizontal cylinder, then the number of square tiles is $n = wh$.
We call \defn{$n$-square-tiled torus} a square-tiled torus tiled by $n$ squares.

Moreover, the lower boundary of the cylinder is necessarily identified with the upper boundary by a translation that respects the square tiling.
Therefore, in the natural coordinates for the horizontal cylinder given by $\RR/w\ZZ \times [0,h]$, this means that there is $t \in \ZZ/w\ZZ$ such that $\RR/w\ZZ \times \{0\}$ is identified with $\RR/w\ZZ \times \{h\}$ by the translation $(x,0) \mapsto (x+t,h)$.
The number $t$ is called the \emph{twist parameter}, and we can always take a representative $t \in \{0,1,\dots,w-1\}$.

Thus, any square-tiled tori uniquely determines the three parameters $w$, $h$ and $t$ as before.
Conversely and similarly, such parameters completely define a square-tiled tori.
See \Cref{f:parameters}.
That is, we have proven the following.

\begin{theorem}[Parametrization]\label{t:parametrization}
There is a one-to-one correspondence between square-tiled tori
and parameters in
\[
\params = \{(w,h,t) \in \NN^2\times \ZZ \colon 0 \leq t < w\},
\]
and between $n$-square-tiled tori
and parameters in
\[
\params_n = \{(w,h,t) \in \params \colon wh = n\}.
\pushQED{\qed}\qedhere\popQED
\]
\end{theorem}

\begin{figure}
\begin{subfigure}{\textwidth}
\centering
\def\ww{9}
\def\www{4}
\def\wwww{5}
\def\hh{3}
\def\tt{2}
\def\tzero{90/\ww}
\def\squaresize{1}
\begin{tikzpicture}[scale=.618,fill=lightgray]
\cylinder
\cylinderparams
\end{tikzpicture}
\hfil
\begin{tikzpicture}[scale=.618,fill=lightgray]
\fundom
\fundomparams
\end{tikzpicture}
\caption{Horizontal cylinder and fundamental domain.}
\end{subfigure}
\begin{subfigure}{\textwidth}
\centering
\def\ww{7}
\def\www{3}
\def\wwww{4}
\def\hh{11}
\def\tt{2}
\begin{tikzpicture}[yscale=1/3,xscale=1/2,fill=lightgray]
\cylinder
\cylinderparams
\end{tikzpicture}
\hfil
\begin{tikzpicture}[scale=1/3,fill=lightgray, thick]
\def\hstep{2*360/(2*\hh+1)}
\def\angle{360/\ww}
\def\angledom{\hstep*\hh/2/\ww}
\def\radius{\squaresize*\ww/2}
\def\Radius{\squaresize*\ww/6}
\def\innerRadius{(\radius-\Radius)}
\def\outerRadius{(\radius+\Radius)}

\filldraw[fill=lightgray!33!white,shift={({-\outerRadius*cos(11*\hstep)},{\outerRadius*sin(11*\hstep)})},rotate=-11*\hstep] (0,0) arc (-180:180:{\Radius} and .4);
\foreach\i in {0,...,11}{
\draw[thick,shift={({-\outerRadius*cos(\i*\hstep)},{\outerRadius*sin(\i*\hstep)})},rotate=-\i*\hstep] (0,0) arc ({180*sign(cos(\i*\hstep))}:0:{\Radius} and {.4*abs(cos(\i*\hstep))});
}

\draw plot[samples=100,domain=0*\angledom:2*\angledom,smooth] ({-(\radius - \Radius*sin(\angle*2 - 2/7*\x))*cos(\x)},{(\radius - \Radius*sin(\angle*2 - 2/7*\x))*sin(\x)});
\draw[very thick] plot[samples=100,domain=0*\angledom+7:10*\angledom,smooth] ({-(\radius - \Radius*sin(\angle*3 - 2/7*\x))*cos(\x)},{(\radius - \Radius*sin(\angle*3 - 2/7*\x))*sin(\x)});
\draw plot[samples=100,domain=0*\angledom+5:14*\angledom+2.5,smooth] ({-(\radius - \Radius*sin(-3*\angle - 2/7*\x))*cos(\x)},{(\radius - \Radius*sin(4*\angle - 2/7*\x))*sin(\x)});
\draw plot[samples=100,domain=0*\angledom:14*\angledom+6.5,smooth] ({-(\radius - \Radius*sin(5*\angle - 2/7*\x))*cos(\x)},{(\radius - \Radius*sin(5*\angle - 2/7*\x))*sin(\x)});
\draw plot[samples=100,domain=5*\angledom:14*\angledom+5,smooth] ({-(\radius - \Radius*sin(6*\angle - 2/7*\x))*cos(\x)},{(\radius - \Radius*sin(6*\angle - 2/7*\x))*sin(\x)});
\draw plot[samples=100,domain=13*\angledom:14*\angledom+.5,smooth] ({-(\radius - \Radius*sin(0*\angle - 2/7*\x))*cos(\x)},{(\radius - \Radius*sin(7*\angle - 2/7*\x))*sin(\x)});

\begin{scope}[even odd rule]
\clip[shift={({-\outerRadius*cos(11*\hstep)},{\outerRadius*sin(11*\hstep)})},rotate=-11*\hstep] {({-.2*cos(11*\hstep)},{.2*sin(11*\hstep)})} arc (-180:180:5.4) (0,0) arc (-180:180:{\Radius} and .4);
\fill[opacity=.8]
({-\innerRadius},0) arc (0:-180:{\Radius} and .4) arc (180:{180-11*\hstep}:{\outerRadius}) --
++ ({180-11*\hstep}:{-2*\Radius}) 
arc ({180-11*\hstep}:180:{\innerRadius});
\end{scope}

\draw[very thick] plot[samples=100,domain=0*\angledom-5:14*\angledom-3.5,smooth] ({-(\radius - \Radius*sin(\angle*1 - 2/7*\x))*cos(\x)},{(\radius - \Radius*sin(\angle*1 - 2/7*\x))*sin(\x)});
\draw plot[samples=100,domain=2*\angledom:14*\angledom-7,smooth] ({-(\radius - \Radius*sin(\angle*2 - 2/7*\x))*cos(\x)},{(\radius - \Radius*sin(\angle*2 - 2/7*\x))*sin(\x)});
\draw[very thick] plot[samples=100,domain=10*\angledom:14*\angledom-4.5,smooth] ({-(\radius - \Radius*sin(\angle*3 - 2/7*\x))*cos(\x)},{(\radius - \Radius*sin(\angle*3 - 2/7*\x))*sin(\x)});
\draw plot[samples=100,domain=0*\angledom-3:5*\angledom,smooth] ({-(\radius - \Radius*sin(-1*\angle - 2/7*\x))*cos(\x)},{(\radius - \Radius*sin(6*\angle - 2/7*\x))*sin(\x)});
\draw plot[samples=100,domain=0*\angledom-7:13*\angledom,smooth] ({-(\radius - \Radius*sin(0*\angle - 2/7*\x))*cos(\x)},{(\radius - \Radius*sin(7*\angle - 2/7*\x))*sin(\x)});

\draw ({-(\radius - \Radius*sin(\angle*1 - 2/7*(-5)))*cos(-5)},{(\radius - \Radius*sin(\angle*1 - 2/7*(-5)))*sin(-5)}) node {$\bullet$};
\draw ({-(\radius - \Radius*sin(\angle*3 - 2/7*(14*\angledom-4.5)))*cos(14*\angledom-4.5)},{(\radius - \Radius*sin(\angle*3 - 2/7*(14*\angledom-4.5)))*sin(14*\angledom-4.5)}) node {$\bullet$};

\draw ({-\innerRadius},0) arc (0:-180:{\Radius} and .4) arc (180:{180-11*\hstep}:{\outerRadius}) ++ ({180-11*\hstep}:{-2*\Radius}) arc ({180-11*\hstep}:180:{\innerRadius});
\foreach\i in {0,...,11}{
\draw[thick,shift={({-\outerRadius*cos(\i*\hstep)},{\outerRadius*sin(\i*\hstep)})},rotate=-\i*\hstep] (0,0) arc ({-180*sign(cos(\i*\hstep))}:0:{\Radius} and {.4*abs(cos(\i*\hstep))});
}
\draw ({-\outerRadius},0) arc (-180:0:{\Radius} and .4);
\draw (0,-7.5) node {};
\end{tikzpicture}
\caption{A square-tiled torus glued from a \emph{twisted} square-tiled cylinder.}
\end{subfigure}
\caption{Parametrization of square-tiled tori.}
\label{f:parameters}
\end{figure}

\begin{remark}
\label{r:sublattices}
From algebraic topology, there is a natural correspondence between square-tiled tori and rank two sublattices of $\ZZ^2 \cong \ZZ[i]$.
In particular, the parametrization in \Cref{t:parametrization} also corresponds to a \emph{canonical} basis for rank two sublattices of $\ZZ^2$ discussed in \Cref{s:parametrization}.
\end{remark}

\subsection*{\texorpdfstring{$\slz$}{SL(2,Z)}-action}
There is a natural action of $\slz$ on square-tiled surfaces, coming from the linear action of $\slz$ on $\RR^2 \cong \CC$.
That is, given $g \in \slz$, we apply $g$ to each of the squares tiling a square-tiled surface $S$ keeping the edge identifications.
This construction gives a new translation surface $gS$ that covers $g(\RR^2)/g(\ZZ^2) = \RR^2/g(\ZZ^2)$.
Since the linear action of $\slz$ on $\RR^2 \cong \CC$ preserves the lattice $\ZZ^2 \cong \ZZ[i]$, we get that the new surface covers $\RR^2/\ZZ^2 \cong \TT$ and therefore $gS$ also is a square-tiled surface.

Note that the linear action of $\slz$ on $\RR^2$ is area preserving.
In particular, the number of squares tiles of a square-tiled surface, which coincides with its area, is preserved by the $\slz$-action.
For example, it is clear from the previous discussion that $\slz$ preserves the square torus $\TT$, the only square-tiled torus that can be tiled with only one unit square.

In this work, we classify $\slz$-orbits using the parametrization of \Cref{t:parametrization}.
In fact, there are two numerical invariants that determine the $\slz$-orbits:
(a) the number $n = wh$ of square tiles, and (b) the number $r = \gcd(w,h,t)$, that corresponds to the side length of the largest square that can tile the corresponding torus, as in \Cref{f:bigger-square-tiles}.
From its definition, and also from its geometrical interpretation, it is clear that $r^2 \divides n$.
We call \defn{$(n,r)$-square-tiled torus} a square-tiled torus with the numerical invariants $n = wh$ and $r = \gcd(w,h,t)$.
We also say that $(n,1)$-square-tiled tori are \defn{cyclic}.

\begin{figure}
\begin{tikzpicture}[scale=1/3,fill=lightgray]
\def\ww{12}
\def\www{3}
\def\wwww{4}
\def\hh{6}
\def\tt{3}
\def\tzero{-3*90/6-90/12}
\def\squaresize{1}
\cylinder
\def\ww{4}
\def\www{2}
\def\wwww{3}
\def\hh{2}
\def\tt{1}
\def\tzero{3*90/12+90/6}
\def\squaresize{3}
\cylinder
\end{tikzpicture}
\hfil
\begin{tikzpicture}[scale=1/3,fill=lightgray]
\def\ww{12}
\def\hh{6}
\def\tt{3}
\def\squaresize{1}
\fundom
\def\ww{4}
\def\hh{2}
\def\tt{1}
\def\squaresize{3}
\fundom
\end{tikzpicture}
\caption{A square-tiled tori tiled with bigger squares.}
\label{f:bigger-square-tiles}
\end{figure}

\begin{remark}
\label{r:cyclic-tori}
The term \emph{cyclic} comes from the fact that the Deck transformations group of the corresponding cover over $\TT$ is a cyclic group.
In \Cref{p:cyclic}, we prove that indeed this is equivalent to $r=1$.
\end{remark}

\begin{theorem}[Orbit classification]\label{t:classification}
The $\slz$-orbits coincide with the sets of $(n,r)$-square-tiled tori, where $n,r \in \NN$ and $r^2 \divides n$.
Moreover, the set of $(n,r)$-square-tiled tori are in one-to-one correspondence with
parameters in
\[
\params_{n,r} = \{(w,h,t) \in \params_n\colon \gcd(w,h,t) = r\}.
\]
\end{theorem}

As a direct consequence we have the following.

\begin{corollary}
Cyclic $n$-square-tiled tori form a single $\slz$-orbit, for every $n\in\NN$.
\qed
\end{corollary}

In particular, there is a unique $\slz$-orbit of $n$-square-tiled tori if and only if every $n$-square-tiled torus is cyclic, and this occurs if and only if $n$ is square-free.

\subsection*{Counting square-tiled tori}
Using \Cref{t:parametrization} and \Cref{t:classification}, we are able to count the total number of $n$-square-tiled tori and the number of square-tiled tori in each $\slz$-orbit.
This can be summarized as follows.

\begin{theorem}[Counting]\label{t:counting}
Let $n \in \NN$.
Then, the total number of $n$-square-tiled tori is
\[
\sigma(n)=\sum_{d\divides n}d
,
\]
where the sum is over \emph{all} divisors of $n$.
Moreover, if $n = r^2 m$ for some $m,r \in \NN$, 
then, the set of $(n,r)$-square-tiled tori form a single $\slz$-orbit of size
\[
\psi(m) = m\cdot\prod_{p\divides m}\left(1+\frac{1}{p}\right)
,
\]
where the product is over the \emph{prime} divisors of $m$.
In particular, the number of cyclic $n$-square-tiled tori is $\psi(n)$.
\end{theorem}

The arithmetic functions $\sigma$ and $\psi$ in \Cref{t:counting} are known as the \defn{sum-of-divisors sigma function} and the \defn{Dedekind psi function}, respectively.
As a side result, we deduce that $\sigma(n)=\sum_{r^2\divides n}\psi(n/r^2)$, where the sum is over \emph{square} divisors of $n$.
We also study the asymptotic behavior of these functions and we conclude the following as a corollary of \Cref{t:counting} and properties of $\psi$ and $\sigma$ (see \Cref{s:asymptotics} for more details and definitions).

\begin{theorem}[Asymptotics]\label{t:asymptotics}
The number $\psi(n)$, of cyclic $n$-square-tiled tori, has the following extremal asymptotic behavior with respect to $\sigma(n)$, the total number of $n$-square-tiled tori:
\[
\liminf_{n\to\infty}\dfrac{\psi(n)}{\sigma(n)} = \dfrac{1}{\zeta(2)}
\qquad\text{and}\qquad
\limsup_{n\to\infty}\dfrac{\psi(n)}{\sigma(n)} = 1
.
\]
Also, the average order of the number of cyclic $n$-square-tiled tori is $1/\zeta(4)$ times the average order of the total number of $n$-square-tiled tori, that is,
\[
\lim_{N\to\infty} \frac{\sum_{n=1}^{N} \psi(n)}{\sum_{n=1}^{N} \sigma(n)} = \frac{1}{\zeta(4)}
.
\]
\end{theorem}

Note a numerical coincidence in relation to deviations from the mean order here and in \cite[Proposition~1.5]{Hubert-Lelievre}: 
\textit{the number of $n$-square-tiled surfaces in $\H(2)$ for prime $n$ is asymptotically $1/\zeta(4)$ times the average order of the number of $n$-square-tiled surfaces in $\H(2)$}.
Here, $\H(2)$ denotes the set of (genus two) translation surfaces with only one conical singularity of cone angle $6 \pi$, that correspond to abelian differentials with only one zero of order \emph{two}.
Our work is in some sense a toy analogue of the work by Hubert--Lelièvre~\cite{Hubert-Lelievre} on square-tiled surfaces in $\H(2)$, which is the best understood moduli space of abelian differentials in higher genus (c.f.~\cite{McMullen:classification-H2}).

% -------------------------------------------------------------------------
\section{Square-tiled tori}

The torus $\TT = \CC/\ZZ[i]$, as any Riemann surface, is a connected, arc-connected and locally simply connected topological space and, in particular, any (connected unbranched) cover over $\TT$ is the quotient of its universal cover $\CC$ by a subgroup of the fundamental group $\pi_1(\TT) \cong \ZZ^2$.
As the fundamental group is abelian, different subgroups determine different (non-isomorphic) covers.
That is, up to isomorphism, each cover of the torus is uniquely determined by a subgroup of $\ZZ^2$.

In particular, each $n$-square-tiled torus is uniquely determined by a sublattice of $\ZZ^2$ of index~$n$.
The cyclic square-tiled tori correspond to primitive sublattices, that is, such that the quotient group is cyclic.

\subsection{Invariants}
Let $\Lambda$ be a (rank two) sublattice of $\ZZ^2$ and $\mathbf{u},\mathbf{v}\in\ZZ^2$ be generators of $\Lambda$.
Note that the number $n = \left|\det(\mathbf{u},\mathbf{v})\right|$ satisfies $n=\left[\ZZ^2:\Lambda\right]$ and, in particular, does not depend on the choice of generators.
Besides, we consider the number $r = \gcd(\mathbf{u},\mathbf{v})$, which also turns out to be independent of the choice of $\mathbf{u}$ and $\mathbf{v}$.
In fact, we have the following.

\begin{lemma} \label{l:1} Let $\Lambda_i = \left<\mathbf{u}_i,\mathbf{v}_i\right>_\ZZ$ be a rank two sublattice of $\ZZ^2$, for $i=1,2$.
Assume that $\Lambda_1\subset\Lambda_2$.
Then, $\gcd(\mathbf{u}_2,\mathbf{v}_2)\divides\gcd(\mathbf{u}_1,\mathbf{v}_1)$.
\end{lemma}
\begin{proof} Since $\Lambda_1\subset\Lambda_2$, there exists $p,q,r,s\in\ZZ$ such that $u_{1j} = p u_{2j} + q v_{2j}$ and $v_{1j} = r u_{2j} + s v_{2j}$.
But $\gcd(\mathbf{u}_2,\mathbf{v}_2)\divides\gcd(u_{2j},v_{2j})$, $\gcd(u_{2j},v_{2j}) \divides p u_{2j} + q v_{2j} = u_{1j}$ and $\gcd(u_{2j},v_{2j}) \divides r u_{2j} + s v_{2j} = v_{1j}$, for $j=1,2$.
We conclude that $\gcd(\mathbf{u}_2,\mathbf{v}_2)\divides\gcd(\mathbf{u}_1,\mathbf{v}_1)$.
\end{proof}

Then $r=\gcd(\mathbf{u},\mathbf{v})$ depends only on $\Lambda$.
Given $\Lambda$, a (rank two) sublattice of $\ZZ^2$, we denote $S(\Lambda)$ the corresponding square-tiled tori $S(\Lambda) \cong \RR^2/\Lambda$.
Similarly, given a square-tiled tori, we denote $\Lambda(S)$ the corresponding sublattice.
Finally, if $S = S(\Lambda)$, we denote the numerical invariants $n(S) = n(\Lambda)$ and $r(S) = r(\Lambda)$ and we say that $S$ is an \defn{$(n,r)$-square-tiled torus}.

\begin{remark} \label{rema:square}
If $S$ is a $(n,r)$-square-tiled surface, then $r^2$ divides $n$.
In fact, if $\mathbf{u},\mathbf{v}$ are generators of $\Lambda=\Lambda(S)$, then $\mathbf{u}_r=\frac{1}{r}\mathbf{u},\mathbf{v}_r=\frac{1}{r}\mathbf{v}\in \ZZ^2$.
In particular, $n/r^2=\left|\det(\mathbf{u},\mathbf{v})\right|/r^2=\left|\det(\mathbf{u}_r,\mathbf{v}_r)\right|$ is integer and $r^2 \divides n$.
\end{remark}

As anticipated by \Cref{t:classification} above, the invariant $r$ plays a key role in this work.
In particular, it allows us to characterize cyclic square-tiled tori.

\begin{proposition}\label{p:cyclic}
Let $S$ be a square-tiled tori.
Then $S$ is cyclic if and only if $r(S) = 1$.
\end{proposition}
\begin{proof}
Let $\Lambda = \Lambda(S)$ and $r=r(S)$.
Let $\mathbf{u},\mathbf{v} \in \ZZ^2$ be generators of $\Lambda$.
Then, $r = \gcd(\mathbf{u},\mathbf{v})$ and we can write $\mathbf{u} = r_ur\mathbf{u}_0$, $\mathbf{v} = r_vr\mathbf{v}_0$, where $r_u = \gcd(\mathbf{u})/r$ and $r_v = \gcd(\mathbf{v})/r$.
In particular, $\gcd(r_u,r_v)= 1$ and $\mathbf{u}_0,\mathbf{v}_0$ are generators of $\ZZ^2$.
Then, since the deck transformations group of the cover $S \to \TT$ is isomorphic to $\ZZ^2/\Lambda$, it is also isomorphic to $\ZZ/r_ur\ZZ \oplus \ZZ/r_vr\ZZ$.
Hence $S$ is cyclic if and only if this last group is cyclic.
But this happens if and only if $\gcd(r_ur,r_vr) = 1$, and $\gcd(r_ur,r_vr) = r\gcd(r_u,r_v) = r$.
\end{proof}

Moreover, for certain values of $n$, there is only cyclic square-tiled tori.
\begin{corollary}
If $n\in\NN$ is square-free, then every $n$-square-tiled tori is cyclic.
\end{corollary}
\begin{proof}
Let $S$ be an $(n,r)$-square-tiled torus.
By \Cref{rema:square}, $r^2$ divides $n$.
But, $n$ is square-free, then $r=1$ and, by \Cref{p:cyclic}, $S$ is cyclic.
\end{proof}

\subsection{Parametrization}\label{s:parametrization}
In this section we give an alternative proof of \Cref{t:parametrization} by means of rank two sublattices of $\ZZ^2$.
In fact, the parameters $w,h,t \in \NN$ for a square-tiled tori determine a \emph{canonical} basis of the corresponding rank two lattices (c.f.
\cite[Section~VII.5.2]{Serre}).

\begin{lemma}
\label{l:parameters}
Let $\Lambda$ be a rank two sublattice of $\ZZ^2$ of index $n\in\NN$.
Then, there exists $w,h\in\NN$ and $t \in \{0,\dots, w-1\}$ such that $\mathbf{u}=(w,0)$ and $\mathbf{v}=(t,h)$ are generators of $\Lambda$.
In particular $n=wh$.
Moreover, these numbers are uniquely determined by $\Lambda$.
\end{lemma}

Recall that the numbers $w, h, t$ have a natural geometric interpretation:
they correspond to natural parameters of the maximal horizontal cylinder filling the square-tiled torus.
The cylinder is isometric to $\RR/w\ZZ \times [0, h]$ and the additional \emph{twist parameter} $t$ carries the information on how to identify the boundaries of the cylinder to recover the original square-tiled torus%
\footnote{A way to agree that we can take $t \in \{0,1,\dots,w-1\}$ is that making an appropriate Dehn twist along the central curve of the cylinder, we can reduce the value of the twist parameter modulo $w$ to obtain diffeomorphic square-tiled tori.}%
.
Namely, the lower boundary $\RR/w\ZZ \times \{0\}$ is identified with the upper boundary $\RR/w\ZZ \times \{h\}$ by the translation $(x,0) \mapsto (x+t,h)$.
See \Cref{f:parameters}.

\begin{proof}[Proof of \Cref{l:parameters}]
Let $\Lambda=\Lambda(S)$ and consider
\begin{align*}
w & = \min\{i\in\NN\colon (i,0)\in\Lambda\},\\
h & = \min\{j\in\NN\colon\exists i\in\ZZ\colon (i,j)\in\Lambda\}, \qquad\text{and} \\
t & = \min\{i\in\NN\colon (i,h)\in\Lambda\}
.
\end{align*}
They are well defined (are finite) since $\Lambda$ is a rank two sublattice of $\ZZ^2$.
Let $\mathbf{u}=(w,0)$ and $\mathbf{v}=(t,h)$.
By definition, $\mathbf{u},\mathbf{v}\in\Lambda$.

We shall prove that $\Lambda\subset\left\langle\mathbf{u},\mathbf{v}\right\rangle_\ZZ$.
Let $(i,j)\in\Lambda$.
Then $j = ch$ for some $c\in\ZZ$.
In fact, if $j = ch+d$ with $c \in \ZZ$ and $d\in\{0,\dots,h-1\}$, then $(i,j)-c(t,h)=(i-ct,d)\in\Lambda$.
By the choice of $h$, we deduce that $d = 0$.
Now, let $i - ct = ew+f$, with $e\in\ZZ$ and $f\in\{0,\dots,w-1\}$.
Then, $(i,j) - c(t,h) - e(w,0) = (f,0)\in\Lambda$ and $w$ is the least positive integer such that $(w,0)\in\Lambda$.
It follows that $f=0$.
Thus, $(i,j)=e\mathbf{u} + c\mathbf{v}$ and $\Lambda\subset\left\langle\mathbf{u},\mathbf{v}\right\rangle_\ZZ$.

Suppose now that, for $\varepsilon \in \{0,1\}$, $\Lambda=\left\langle\mathbf{u}_\varepsilon,\mathbf{v}_\varepsilon\right\rangle_\ZZ$ with $\mathbf{u}_\varepsilon=(w_\varepsilon,0)$ and $\mathbf{v}_\varepsilon=(t_\varepsilon,h_\varepsilon)$.
Then, there exist $p_\varepsilon,q_\varepsilon,r_\varepsilon,s_\varepsilon\in\ZZ$ such that $\mathbf{u}_\varepsilon=p_\varepsilon\mathbf{u}_{1-\varepsilon}+q_\varepsilon\mathbf{v}_{1-\varepsilon}$ and $\mathbf{v}_\varepsilon=r_\varepsilon\mathbf{u}_{1-\varepsilon}+q_\varepsilon\mathbf{v}_{1-\varepsilon}$.
That is, $w_\varepsilon = p_\varepsilon w_{1-\varepsilon} + q_\varepsilon t_{1-\varepsilon}$, $0 = q_\varepsilon h_{1-\varepsilon}$, $t_\varepsilon = r_\varepsilon w_{1-\varepsilon} + r_\varepsilon t_{1-\varepsilon}$ and $h_\varepsilon = s_\varepsilon h_{1-\varepsilon}$.
But then, $q_\varepsilon=0$, $p_\varepsilon=1$, $s_\varepsilon=1$ and $r_\varepsilon=0$, and we conclude that $(w_0,h_0,t_0)=(w_1,h_1,t_1)$.
\end{proof}

\begin{proof}[Proof of \Cref{t:parametrization}]
Let $n \in \NN$.
Then, $n$-square-tiled tori are in one-to-one correspondence with rank two sublattices of $\ZZ^2$ of index $n$ and, by \Cref{l:parameters}, these are in one-to-one correspondence with parameters $w,h \in \NN$ and $t \in \{0,\dots, w-1\}$, such that $wh = n$.
\end{proof}

In what follows, we denote by $S(w,h,t)$ the square-tiled torus associated with $(w,h,t) \in \params$.

\section{\texorpdfstring{$\slz$}{SL(2,Z)}-action}\label{s:SL2Z-action}

In this section we describe the $\slz$-action on square-tiled tori and classify their orbits, proving \Cref{t:classification}.

Recall that $\slz$ acts as the linear action in $\CC \cong \RR^2$.
In terms of lattices, we have that if $A\in\slz$ and $S$ is a square-tiled torus, then $\Lambda(A\cdot S) = A\cdot\Lambda(S)$.
This action clearly preserves the area of square-tiled tori, that is, the number of squares.
In particular, $\slz$-orbits are finite (for the counting see \Cref{s:counting}).

In order to prove \Cref{t:classification}, we use standard generators of $\slz$.
Let $R = \begin{psmallmatrix}0&-1\\1&0\end{psmallmatrix}$ and $U = \begin{psmallmatrix}1&1\\0&1\end{psmallmatrix}$.
One can see, using row and column operations and the Euclidean algorithm, that $\slz$ is generated by $R$ and $U$.

\begin{lemma}\label{l:r-invariant}
The number $r$ is an $\slz$-invariant, that is, $r(S)=r(A\cdot S)$, for any square-tiled torus $S$ and $A\in\slz$.
In particular, being a cyclic square-tiled torus is an $\slz$-invariant property.
\end{lemma}
\begin{proof}
This follows from the $\slz$-invariance of $\gcd$.
In fact, for $u,v\in\ZZ$, we have that $\gcd(u,v)=\gcd(u,u+v)$ and $\gcd(u,v)=\gcd(v,-u)$.
Then $\gcd$ is $U$ and $R$-invariant and hence, $\slz$-invariant.

Now, let $S$ be a square-tiled surface and $\Lambda = \Lambda(S)$.
If $\mathbf{u},\mathbf{v} \in \ZZ^2$ are generators of $\Lambda$, then $A\cdot\mathbf{u}$ and $A\cdot\mathbf{v}$ are generators of $A\cdot \Lambda = \Lambda(A\cdot S)$.
Then, by the $\slz$-invariance of $\gcd$, we get that $\gcd(A\cdot\mathbf{u},A\cdot\mathbf{v}) = \gcd(\mathbf{u},\mathbf{v})$, that is, $r(A\cdot S) = r(S)$.

The claim about cyclic square-tiled tori follows then from \Cref{p:cyclic}.
\end{proof}

In geometric terms, the action of $R$ and $U$ can be understood by means of their action on a square in $\RR^2$, as in \Cref{f:action}.
Then, the action on square-tiled surfaces is obtained by applying this to all the square tiles and keeping the corresponding side identifications.
As an illustration of this, in \Cref{f:SL2Z-orbit} we exhibit the $\slz$-orbit of the cyclic $4$-square-tiled tori using the action of the generators $R$ and $U$ of $\slz$.

\begin{figure}
\begin{tikzpicture}[scale=1,fill=lightgray,thick]
\begin{scope}[shift={(-1.5,0)}]
\filldraw (0,0) rectangle (1,1);
\draw (0,0) node {$\bullet$};
\draw[very thick, dots] (-.25,0) -- (1.25,0) (-.25,1) -- (1.25,1) (0,-.25) -- (0,1.25) (1,-.25) -- (1,1.25);
\draw[->] (.6,0) -- (.618,0);
\draw[->>] (0,.6) -- (0,.618);
\end{scope}
\draw[->] (0,.5) -- node[midway,above] {$R$} (1,.5);
\begin{scope}[shift={(1.5,0)}]
\filldraw (0,0) rectangle (1,1);
\draw (1,0) node {$\bullet$};
\draw[very thick, dots] (-.25,0) -- (1.25,0) (-.25,1) -- (1.25,1) (0,-.25) -- (0,1.25) (1,-.25) -- (1,1.25);
\draw[->>] (1-.6,0) -- (1-.618,0);
\draw[->] (1,.6) -- (1,.618);
\end{scope}
\end{tikzpicture}
\hfil
\begin{tikzpicture}[scale=1,fill=lightgray,thick]
\begin{scope}[shift={(-1.5,0)}]
\filldraw (0,0) rectangle (1,1);
\draw (0,0) node {$\bullet$};
\draw[very thick, dots] (-.25,0) -- (1.25,0) (-.25,1) -- (1.25,1) (0,-.25) -- (0,1.25) (1,-.25) -- (1,1.25);
\draw[->] (.6,0) -- (.618,0);
\draw[->>] (0,.6) -- (0,.618);
\end{scope}
\draw[->] (0,.5) -- node[midway,above] {$U$} (1,.5);
\begin{scope}[shift={(1.5,0)}]
\filldraw (0,0) -- (1,0) -- (2,1) -- (1,1) -- cycle;
\draw (0,0) node {$\bullet$};
\draw[very thick, dots] (-.25,0) -- (2.25,0) (-.25,1) -- (2.25,1) (0,-.25) -- (0,1.25) (1,-.25) -- (1,1.25) (2,-.25) -- (2,1.25);
\draw[->] (.6,0) -- (.618,0);
\draw[->>] (0+.6,.6) -- (0+.618,.618);
\end{scope}
\end{tikzpicture}
\caption{Action of the standard generators of $\slz$ on squares.}
\label{f:action}
\end{figure}

\begin{figure}
\begin{tikzpicture}[scale=.4,fill=lightgray,thick]
\def\ww{1}
\def\hh{4}
\def\tt{0}
\begin{scope}[shift={(0,-{\hh/2})}]
\draw[->] (-.5,1.5) arc (280:80:.5) node[midway,left] {$U$};
\fundom
\draw[<->] (1.5,2) -- node[midway,above] {$R$} (2.5,2);
\end{scope}
\def\ww{4}
\def\hh{1}
\def\tt{0}
\begin{scope}[shift={(3,-{\hh/2})}]
\fundom
\draw[->] (4,1) ++(45:.5) -- node[midway,above left] {$U$} ++(45:1);
\draw[<-] (4,0) ++(-45:.5)  -- node[midway,below left] {$U$} ++(-45:1);
\end{scope}

\begin{scope}[shift={({10+sqrt(2)},{-\hh/2})}]
\draw[<->] (0,0) -- node [midway,right] {$R$} (0,1);
\def\tt{1}
\begin{scope}[shift={({-(\ww+\tt)/2},{sqrt(2)+\hh})}]
\fundom
\end{scope}
\def\tt{3}
\begin{scope}[shift={({-(\ww+\tt)/2},{-sqrt(2)-\hh})}]
\fundom
\end{scope}
\end{scope}

\def\tt{2}
\begin{scope}[shift={({13+2*sqrt(2)},-{\hh/2})}]
\fundom
\draw[<-] (0,1) ++(135:.5) -- node[midway,above right] {$U$} ++(135:1);
\draw[->] (0,0) ++(-135:.5)  -- node[midway,below right] {$U$} ++(-135:1);
\end{scope}

\def\ww{2}
\def\hh{2}
\def\tt{1}
\begin{scope}[shift={({21+2*sqrt(2)},-{\hh/2})}]
\draw[->] (3.5,.5) arc (-100:100:.5) node[midway,right] {$U$};
\fundom
\draw[<->] (-.5,1) -- node[midway,above] {$R$} (-1.5,1);
\end{scope}
\end{tikzpicture}
\caption{The $\slz$-orbit of cyclic $4$-square-tiled tori.}
\label{f:SL2Z-orbit}
\end{figure}

Now, let $S$ be a square-tiled torus with parameters $(w,h,t) \in \params$, that is, $S = S(w,h,t)$.
\begin{itemize}
\item
The action of  $R$ corresponds to a rotation by $\pi/2$ (anticlockwise).
Then, the parameters of $R\cdot S$ correspond to the \emph{vertical cylinder}%
\footnote{Similarly to the description of the \emph{horizontal cylinder} filling the square-tiled torus and giving rise to the parameters $w$, $h$ and $t$, one can consider the cylinder decomposition in the vertical direction.
Recall that, in genus one, any cylinder decomposition consists in only one cylinder with its boundary identified by translation.}
of $S$.
The width and height parameters are $w_R = hw/\gcd(w,t)$ and $h_R = \gcd(w,t)$ respectively.
The twist parameter $t_R$ of the vertical cylinder decomposition is also computable, but needless for our purposes.
However, by \Cref{l:r-invariant}, 
$r$ is an $\slz$-invariant and, in particular 
$r$ divides all parameters of square-tiled tori in the same orbit.
In particular, $r$ divides the twist parameter $t_R$ (and this is the only fact we need to know about $t_R$).
\item
On the other hand, $U$ preserves the horizontal direction.
Its action preserves the width $w$ and the height $h$, and changes the twist parameter $t$ to $t_U=(t + h) \bmod w$.
In particular, the size of the $U$-orbit of $S$ is $w/\gcd(w,h)$ and we say that the \defn{canonical representative} of the $U$-orbit of $S$ is the square-tiled torus $S(w,h,t \bmod\gcd(w,h))$.
\end{itemize}
In summary, the $\slz$-action in the parameter space $\params$ is completely determined by the action of $R$ and $U$ which correspond to $R\cdot(w,h,t)=(w_R = hw/\gcd(w,t),h_R = \gcd(w,t),t_R)$ and $U\cdot(w,h,t)=(w,h,t_U=(t + h) \bmod w)$.

\begin{lemma}\label{l:r-complete}
Let $n,r\in\NN$ with $r^2 \divides n$ and assume that $S_1,S_2$ are two $(n,r)$-square-tiled tori.
Then, $S_1$ and $S_2$ are in the same $\slz$-orbit, that is, there exist $A\in\slz$ such that $S_2=A\cdot S_1$.
\end{lemma}
\begin{proof}
Since $r^2 \divides n$, $(n/r,r,0) \in \params_n$.
Let $S_{n,r} = S(n/r,r,0)$.
Then, it is enough to show that $S_{n,r}$  is in the $\slz$-orbit of any $(n,r)$-square-tiled torus $S$.

For this, let $(w,h,t) \in \params$ such that $S = S(w,h,t)$.
In particular, $r = \gcd(w,h,t)$.
\begin{enumcases}
\item \label{s:case-1}
If $w=n/r$ (and $h=r$), we can apply $U^k$ to $S$, with $k = (n-tr)/r^2$, to get $S_{n,r} = U^k S$.
In particular, $S_{n,r} \in \slz\cdot S$.
Note that $k$ is integer since $r^2 \divides n$ and $r \divides t$.
\item
If $w\neq n/r$, let $S_U = S(w,h,t_U)$ be the canonical representative of the $U$-orbit of $S$, that is, $t_U=t \bmod \gcd(w,h)$.
Then, $S_U \in \slz \cdot S$ as well.
\begin{enumcases}
\item \label{s:case-2a}
If $t_U = 0$, then $r=\gcd(w,h)$ and $S(w,h,r)$ is in the $U$-orbit of $S$.
We apply $R$ to $S(w,h,r)$ to obtain $S(wh/r,r,*) \in \slz\cdot S$.
Then, by \cref{s:case-1} applied to $S(wh/r,r,*)$, we get $S_{n,r}\in \slz \cdot S$.
\item \label{s:case-2b}
If $t_U \neq 0$, let $S' = S(w',h',t')$ be the canonical representative of the $U$-orbit of $R\cdot S_U$, that is, $w' = hw/\gcd(w,t_U)$, $h' = \gcd(w,t_U)$ and $t' < \gcd(w',h')$.
\end{enumcases}
\end{enumcases}

Then, either we are in \cref{s:case-1} or \cref{s:case-2a} and $S_{n,r} \in \slz\cdot S$, or we are in \cref{s:case-2b} and there is $S' = S(w',h',t') \in \slz\cdot S$ with 
\[
h' = \gcd(w,t_U) \leq t_U < h
\qquad\text{and}\qquad
w' = h\frac{w}{\gcd(w,t_U)}=w\frac{h}{\gcd(w,t_U)} \geq w\frac{h}{t_U} > w
.
\]

That is, either $S_{n,r} \in \slz\cdot S$ or there is $S' = S(w',h',t') \in \slz\cdot S$ with $w' > w$.
It follows that, if $S_{n,r} \notin \slz\cdot S$, we can produce a sequence $(S_k = S(w_k,h_k,t_k))_{k\in\NN}$ in $\slz\cdot S$ with strictly increasing width parameters: $w_{k+1} > w_k$, $k \in \NN$.
But this is impossible as $\slz$-orbits are finite (and this can be seen as, for example, $n = wh$ is an $\slz$-invariant).
Then, $S_{n,r}$ is in the $\slz$-orbit of $S$ and we conclude that all $(n,r)$-square-tiled tori are in the same $\slz$-orbit.
\end{proof}

\begin{proof}[Proof of \Cref{t:classification}]
Let $S_1$ and $S_2$ be two square-tiled tori.
If $S_1$ and $S_2$ are in the same $\slz$-orbit, then $n(S_1) = n(S_2)$ and, by \Cref{l:r-invariant}, $r(S_1) = r(S_2)$.
Conversely, if $n(S_1) = n(S_2)$ and $r(S_1) = r(S_2)$ then, by \Cref{l:r-complete}, $S_1$ and $S_2$ are in the same $\slz$-orbit.
\end{proof}

\section{Counting square-tiled tori}\label{s:counting}

In this section we count the total number of $n$-square-tiled tori as well as the number of square-tiled tori in each $\slz$-orbit, proving \Cref{t:counting}.
For this purpose, we use the correspondence between the set of $n$-square-tiled tori and $\params_n$ given by \Cref{t:parametrization} and the correspondence between $\slz$-orbits and $\params_{n,r}$, given by \Cref{t:classification}.
In particular, we show that the size $\psi(n)$ of the only $\slz$-orbit of cyclic $n$-square-tiled tori depends on the prime factors of $n$.

The following corresponds to the first statement of \Cref{t:counting}.

\begin{proposition}\label{p:counting-all}
Let $n \in \NN$.
Then, the total number of $n$-square-tiled tori is
\[
\sigma(n)=\sum_{d\divides n}d
,
\]
where the sum is over \emph{all} divisors of $n$.
\end{proposition}
\begin{proof}
By \Cref{t:parametrization}, there is a one-to-one correspondence between the set of $n$-square-tiled tori and $\params_n = \{(w,h,t) \in \NN^2\times\ZZ\colon wh = n, 0 \leq t < w\}$.
Then, the number of $n$-square-tiled tori is
\[
\# \params_n = \sum_{wh = n} w = \sum_{w\divides n} w = \sigma(n).
\]

\end{proof}

In order to prove the second statement of \Cref{t:counting}, we use some basics of number theory and, more precisely, of \defn{arithmetic functions}, that is, functions from $\NN$ to $\RR$.
An arithmetic function $\Phi:\NN\to\RR$ is said to be \defn{multiplicative} if, whenever $n_1,n_2\in\NN$ are coprimes, that is, $\gcd(n_1,n_2)=1$, then $\Phi(n_1n_2)=\Phi(n_1)\Phi(n_2)$.
In particular, we use the arithmetic function
\[
\varphi(n) = \#\!\left\{k\in\{0,\dots,n-1\}\colon \gcd(n,k)=1\right\}
,
\]
known as the \defn{Euler's totient function}, which is multiplicative and satisfies the \emph{Euler's product formula} (see \cite[Section~5.5]{Hardy-Wright})
\begin{formula}\label{f:Euler}
\varphi(n)=n\prod_{p\divides n}\left(1-\frac{1}{n}\right)
,
\end{formula}
where the product is over the \emph{prime} divisors of $n\in\NN$.

We begin with the case of cyclic square-tiled tori.
\begin{lemma} \label{l:cyclic-euler}
Let $n\in\NN$.
Then, the number of cyclic $n$-square-tiled tori is equal to
\[
\psi(n) = \sum_{wh=n}\frac{w}{\gcd(w,h)}\cdot\varphi(\gcd(w,h)).\]
\end{lemma}
\begin{proof}
By \Cref{t:classification}, there is a one-to-one correspondence between the set of cyclic $n$-square-tiled tori and $\params_{n,1} = \{(w,h,t) \in \params_n\colon \gcd(w,h,t) = 1\}$.
That is, the number of cyclic $n$-square-tiled tori coincides with
\[
\#\params_{n,r} = \sum_{wh=n}\#\!\left\{t\in\{0,\dots,w-1\}\colon \gcd(w,h,t)=1\right\}
.\]

As discussed in \Cref{s:SL2Z-action}, each square-tiled torus $S = S(w,h,t)$ is in the $U$-orbit of a unique canonical representative $S_U = S(w,h,t_U)$ with $0\leq t_U < \gcd(w,h)$, and there are exactly $w/\gcd(w,h)$ square-tiled tori in the corresponding $U$-orbit.
Then, 
\begin{align*}
\#\params_{n,r} 
& = \sum_{wh=n}\frac{w}{\gcd(w,h)}\#\!\left\{t_U\in\{0,\dots,\gcd(w,h)-1\}\colon\gcd(w,h,t_U) = 1\right\}
\\
& = \sum_{wh=n}\frac{w}{\gcd(w,h)}\#\!\left\{t\in\{0,\dots,\gcd(w,h)-1\}\colon\gcd(\gcd(w,h),t) = 1\right\}
\\
& = \sum_{wh=n}\frac{w}{\gcd(w,h)}\varphi(\gcd(w,h))
= \psi(n)
.
\qedhere
\end{align*}
\end{proof}

\begin{lemma} \label{l:mult}
The function $n\mapsto\psi(n)$ is a multiplicative arithmetic function such that, for prime $p$ and positive integer $\alpha$,
\[
\psi(p^\alpha)=p^\alpha\cdot\left(1+\frac{1}{p}\right).\]
\end{lemma}
\begin{proof}
Let $n_1,n_2\in\NN$ be coprimes and $n=n_1n_2$.
We want to show that $\psi(n)=\psi(n_1)\psi(n_2)$.
Let $w,h$ be such that $wh=n$.
As $n_1$ and $n_2$ are coprimes, we can decompose $w=w_1w_2$ and $h=h_1h_2$ where $w_i,h_i \divides n_i$, for $i=1,2$.
Moreover, $w_ih_i=n_i$, for $i=1,2$, and this decomposition is unique.
Furthermore, since $\gcd(n_1,n_2)=1$, then we also have that $\gcd(w_i,w_j) = \gcd(w_i,h_j)=\gcd(h_i,h_j)=1$, for $i\neq j \in \{1,2\}$.
It follows that
\[
\gcd(w,h) = \gcd(w_1w_2,h_1h_2) = \gcd(w_1,h_1h_2)\gcd(w_2,h_1h_2) = \gcd(w_1,h_1)\gcd(w_2,h_2)
.
\]
But then,
\begin{align*}
\psi(n) &= \sum_{wh=n}\frac{w}{\gcd(w,h)}\cdot\varphi(\gcd(w,h)) \\
&= \sum_{w_1h_1=n_1}\sum_{w_2h_2=n_2}\frac{w_1w_2}{\gcd(w_1,h_1)\gcd(w_2,h_2)}\cdot\varphi(\gcd(w_1,h_1)\gcd(w_2,h_2)) \\
&= \sum_{w_1h_1=n_1}\sum_{w_2h_2=n_2}\frac{w_1w_2}{\gcd(w_1,h_1)\gcd(w_2,h_2)}\cdot\varphi(\gcd(w_1,h_1))\varphi(\gcd(w_2,h_2)) \\
&= \sum_{w_1h_1=n_1}\sum_{w_2h_2=n_2}\left(\frac{w_1}{\gcd(w_1,h_1)}\cdot\varphi(\gcd(w_1,h_1))\right)\left(\frac{w_2}{\gcd(w_2,h_2)}\cdot\varphi(\gcd(w_2,h_2))\right) \\
&= \left(\sum_{w_1h_1=n_1}\frac{w_1}{\gcd(w_1,h_1)}\cdot\varphi(\gcd(w_1,h_1))\right) \hspace{-2pt} \left(\sum_{w_2h_2=n_2}\frac{w_2}{\gcd(w_2,h_2)}\cdot\varphi(\gcd(w_2,h_2))\right) \\
&= \psi(n_1)\psi(n_2),
\end{align*}
where we have used the fact that the Euler $\varphi$ function is multiplicative and that $\gcd(w_1,h_1)$ and $\gcd(w_2,h_2)$ are coprimes since $\gcd(w_1,h_1) \divides n_1 = w_1 h_1$ and $\gcd(w_2,h_2) \divides n_2 = w_2 h_2$, and $n_1$, $n_2$ are coprimes by hypothesis.

Let $p$ be prime and $\alpha$ a positive integer.
Since $p$ is the only prime factor of $p^\alpha$, using the Euler's product \cref{f:Euler}, we have
\begin{multline*}
\psi(p^\alpha) = \sum_{wh=p^\alpha}\frac{w}{\gcd(w,h)}\cdot\varphi(\gcd(w,h)) 
= 1 + \sum_{k=1}^{\alpha-1}p^k\left(1-\frac{1}{p}\right) + p^\alpha \\
= 1 + (p-1)\sum_{k=1}^{\alpha-1}p^{k-1} + p^\alpha 
= 1 + (p-1)\sum_{k=0}^{\alpha-2}p^k + p^\alpha 
= 1 + (p-1)\frac{p^{\alpha-1}-1}{p-1} + p^\alpha \\
= p^{\alpha-1} + p^\alpha = p^{\alpha-1}(1+p) 
= p^\alpha\cdot\left(1+\frac{1}{p}\right).
\qedhere
\end{multline*}
\end{proof}

Multiplicative arithmetic functions are completely defined from their values in prime powers.
\begin{corollary}
Let $n \in \NN$.
The number of cyclic $n$-square-tiled tori is
\[
\psi(n) = n\cdot\prod_{p\divides n}\left(1+\frac{1}{p}\right),
\]
where the product is over the distinct \emph{prime} numbers dividing $n$.
\end{corollary}
\begin{proof}
Let $n$ be a positive integer, and $n=\prod_{i=1}^{\omega(n)}p_i^{\alpha_i}$ its prime decomposition.
Then
\[
\psi(n) = \psi\left(\prod_{i=1}^{\omega(n)}p_i^{\alpha_i}\right) 
= \prod_{i=1}^{\omega(n)}\psi(p_i^{\alpha_i}) 
= \prod_{i=1}^{\omega(n)}p_i^{\alpha_i}\left(1+\frac{1}{p_i}\right) \\
= n\cdot\prod_{i=1}^{\omega(n)}\left(1+\frac{1}{p_i}\right).
\]
\end{proof}

In \Cref{s:asymptotics} we exhibit the asymptotic behaviors of the arithmetic functions $\psi$ and $\sigma$.
In particular, proving \Cref{t:asymptotics}.

Analogously to \Cref{l:cyclic-euler}, we now consider the general case of $(n,r)$-square-tiled tori and deduce the following.
\begin{proposition}\label{p:counting-orbits}
Let $n,r \in\NN$.
Then, the number of $(n,r)$-square-tiled tori is $\psi(n/r^2)$ if $r^2\divides n$ and zero, otherwise.
\end{proposition}
\begin{proof} We shall show that the number $\psi_r(n) = \# \params_{n,r}$ of triples $(w,h,t)$ such that $wh=n$, $0\leq t < w$ and $\gcd(w,h,t)=r$, is equal to $\psi(n/r^2)$.
Analogously to the proof of \Cref{l:cyclic-euler},
\[
\psi_r(n) = \sum_{wh=n}\frac{w}{\gcd(w,h)}\#\!\left\{t\in\{0,\dots,\gcd(w,h)-1\}\colon\gcd(w,h,t)=r\right\}.\]
In particular, the summands are non-zero if and only if $r\divides w,h$, and in those cases, we write $w=r\tilde{w}$, $h=r\tilde{h}$.
In particular, if $r^2$ does not divide $n$, there are no such $w,h$ with $wh=n$ and then, $\psi_r(n)=0$.
If $r^2$ divides $n$, then
\begin{align*}
\psi_r(n)
& = \sum_{r\tilde{w}r\tilde{h}=n}\frac{r\tilde{w}}{\gcd(r\tilde{w},r\tilde{h})}\#\!\left\{t\in\{0,\dots,\gcd(r\tilde{w},r\tilde{h})-1\}\colon\gcd(r\tilde{w},r\tilde{h},t)=r\right\} \\
& = \sum_{\tilde{w}\tilde{h}=\frac{n}{r^2}}\frac{r\tilde{w}}{r\gcd(\tilde{w},\tilde{h})}\#\!\left\{r\tilde{t}\in\{0,\dots,\gcd(r\tilde{w},r\tilde{h})-1\}\colon\gcd(r\tilde{w},r\tilde{h},r\tilde{t})=r\right\} \\
& = \sum_{\tilde{w}\tilde{h}=\frac{n}{r^2}}\frac{\tilde{w}}{\gcd(\tilde{w},\tilde{h})}\#\!\left\{\tilde{t}\in\{0,\dots,\gcd(\tilde{w},\tilde{h})-1\}\colon\gcd(\tilde{w},\tilde{h},\tilde{t})=1\right\} \\
& = \psi\left(\frac{n}{r^2}\right).
\qedhere
\end{align*}
\end{proof}

\begin{proof}[Proof of \Cref{t:counting}]
The first statement of \Cref{t:counting} about the total number of $n$-square-tiled surfaces is \Cref{p:counting-all}, and the second, about $(n,r)$-square-tiled tori, follows from \Cref{p:counting-orbits}.
\end{proof}

As side result, we obtain the following number theoretic result.
\begin{corollary}
Let $n\in\NN$.
Then
\[
\sigma(n) = \sum_{r^2\divides n}\psi\left(\frac{n}{r^2}\right),
\]
where the sum is over \emph{square} divisors of~$n$.
\end{corollary}
\begin{proof}
By \Cref{p:counting-all}, we know that $\sigma(n)$ is the size of the set of $n$-square-tiled tori.
Now, an $n$-square-tiled torus is an $(n,r)$-square-tiled torus for some $r\in\NN$ with $r^2\divides n$ and the number of $(n,r)$-square-tiled tori is $\psi(n/r^2)$.
\end{proof}

Furthermore, it is also possible to determine the asymptotics of the size of any $\slz$-orbit, since the pair $(n,r)$ classifies them and we know the asymptotics for $\psi$ (see \Cref{s:asymptotics} and \Cref{t:asymptotics}).

\section{Asymptotics for sum-of-divisors sigma and Dedekind psi functions} \label{s:asymptotics}

In this section we study the asymptotic behavior of the two arithmetic functions relevant to the counting of square-tiled tori, proving \Cref{t:asymptotics}.
Recall that, by \Cref{t:counting}, the number of $n$-square-tiled tori corresponds to the \defn{sum-of-divisors sigma function}
\[
\sigma(n)=\sum_{d\divides n}d
.
\]
On the other hand, the number of cyclic $n$-square-tiled tori is given by the \defn{Dedekind psi function}
\[
\psi(n) = n\cdot\prod_{p\divides n}\left(1+\frac{1}{p}\right)
.
\]

By \defn{arithmetic function} we mean any function $f \colon \NN \to \RR$.
We say that $f$ has \defn{minimal order} $m \colon \NN \to \RR$ and \defn{maximal order} $M \colon \NN \to \RR$ if
\[
\liminf_{n \to \infty} \frac{f(n)}{m(n)} = 1
\qquad\text{and}\qquad
\limsup_{n \to \infty} \frac{f(n)}{M(n)} = 1
,
\]
and are referred to as \defn{extremal orders} of $f$.
We also say that $f$ has \defn{average order} $a \colon \NN \to \RR$ if
\[
\lim_{N \to \infty} \frac{\sum_{n=1}^{N} f(n)}{\sum_{n=1}^{N} a(n)} = 1
.
\]

The asymptotic behavior of sum-of-divisors sigma function is well known since Gr\"onwall's theorem in 1913~\cite{Gronwall}.
Here we present a summary of asymptotics for sum-of-divisors sigma function (see \cite[Section~18.3]{Hardy-Wright}).
\begin{theorem}\label{t:a1}
The sum-of-divisors sigma function satisfy the following asymptotic behavior:
\begin{enumerate}
\item \indent\rlap{Minimal order $n$:}\hfil
$\displaystyle\liminf_{n\to\infty}\frac{\sigma(n)}{n} = 1$.
\item \indent\rlap{Maximal order $e^\gamma n\ln\ln n$:}\hfil
$\displaystyle\limsup_{n\to\infty}\frac{\sigma(n)}{n\ln\ln n} = e^\gamma$.
\item \indent\rlap{Average order $\zeta(2)n$:}\hfill\phantom{\qedsymbol\qedsymbol}
$\displaystyle\lim_{\mathclap{N\to\infty}}\frac{\sum_{n=1}^{N} \sigma(n)}{\sum_{n=1}^{N}n} = \zeta(2)$.%
\pushQED{\qed}\qedhere\popQED
\end{enumerate}
\end{theorem}

The extremal values of Dedekind psi function are studied in~\cite{Sole-Planat} and the mean order in~\cite{Perez}.
In the following we present the extremal and average orders of the Dedekind psi function relative to sum-of-divisors sigma function.
In particular, by \Cref{t:counting}, the following is equivalent to \Cref{t:asymptotics}.

\begin{theorem} \label{t:a2}
The Dedekind psi function has the following extremal asymptotic behavior with respect to sum-of-divisors sigma function:
\[
\liminf_{n\to\infty}\frac{\psi(n)}{\sigma(n)} = \frac{1}{\zeta(2)}
\qquad\text{and}\qquad
\limsup_{n\to\infty}\frac{\psi(n)}{\sigma(n)} = 1
.
\]
Also, the average order of the Dedekind psi function is $1/\zeta(4)$ times the average order of the sum-of-divisors sigma function, that is,
\[
\lim_{N\to\infty} \frac{\sum_{n=1}^{N} \psi(n)}{\sum_{n=1}^{N} \sigma(n)} = \frac{1}{\zeta(4)}
.
\]
\end{theorem}
\begin{proof}
Let $n$ be a positive integer, and $n=\prod_{i=1}^{\omega(n)}p_i^{\alpha_i}$ its prime decomposition.
The sum-of-divisors sigma function satisfies (see \cite[Section~16.7]{Hardy-Wright})
\[
\sigma(n)=\prod_{i=1}^{\omega(n)}\frac{p_i^{\alpha_i+1}-1}{p_i-1}
\]
and recall that Dedekind psi function is given by
\[
\psi(n)=\prod_{i=1}^{\omega(n)}p_i^{\alpha_i}+p_i^{\alpha_i-1}
.
\]
Then, defining $\rho(n) = \psi(n)/\sigma(n)$, we have that
\[
\rho(n)
= \prod_{i=1}^{\omega(n)}\frac{p_i^{\alpha_i+1}+p_i^{\alpha_i-1}}{p_i^{\alpha_i+1}-1} 
= \prod_{i=1}^{\omega(n)}\left(1-\frac{1}{p_i^2}\right)\cdot \prod_{i=1}^{\omega(n)}{\left(1-\frac{1}{p_i^{\alpha_i+1}}\right)}^{-1}.
\]
Since $p_i^{\alpha_i+1}\geq p_i^2$, $\rho(n)\leq 1$ with equality if and only if $\alpha_i=1$ for all $i=1,\dots,\omega(n)$, that is, if and only if $n$ is square-free. In particular, $\limsup \rho(n) = 1$.
Furthermore,
\[
\prod_{i=1}^{\omega(n)}\left(1-\frac{1}{p_i^2}\right)\geq \prod_{p \text{ prime}}\left(1-\frac{1}{p^2}\right)
= \frac{1}{\zeta(2)}
\qquad\text{and}\qquad
\prod_{i=1}^{\omega(n)}\left(1-\frac{1}{p_i^{\alpha_i+1}}\right) \leq 1
.
\]
Therefore, $\rho(n)\geq 1/\zeta(2)$. Taking $n_k:=(p_k\#)^k$ to be the $k$th power of the $k$th primorial (product of the first $k$ primes) we have
\begin{align*}
\lim_{k\to\infty} \rho(n_k) &= \lim_{k\to\infty} \prod_{i=1}^{k}\left(1-\frac{1}{p_i^2}\right)\cdot \prod_{i=1}^{k}{\left(1-\frac{1}{p_i^{k+1}}\right)}^{-1} \\
& = \prod_{p \text{ prime}}\left(1-\frac{1}{p^2}\right) \cdot \lim_{k\to\infty} \prod_{i=1}^{k}{\left(1-\frac{1}{p_i^{k+1}}\right)}^{-1} \\
&\leq  \frac{1}{\zeta(2)}\cdot \lim_{k\to\infty}\prod_{i=1}^{k}{\left(1-\frac{1}{2^{k+1}}\right)}^{-1}
= \frac{1}{\zeta(2)}\cdot \lim_{k\to\infty}{\left(1-\frac{1}{2^{k+1}}\right)}^{-k}
= \frac{1}{\zeta(2)},
\end{align*}
and $\liminf \rho(n) = 1/\zeta(2)$.

For the mean order, let $\q$ be the indicator function of square-free integers and denote
\[
\Psi(n) = n\sum_{d|n}\frac{\q(d)}{d}
,
\]
where the sum is over \emph{all} divisors of $n$.
We claim that $\Psi = \psi$. In fact, is easy to see that $\q$ is multiplicative and analogous computations to those in the proof of \Cref{l:mult} shows that $\Psi$ is multiplicative. But multiplicative functions are completely defined from its values in prime powers and
\[
\Psi(p^\alpha) = p^\alpha \sum_{\beta=1}^{\alpha}\frac{\q(p^\beta)}{p^\beta} = p^\alpha\left(1 + p^{-1}\right) = \psi(p^\alpha),
\]
proving the claim.
We also use the fact (see \cite[Section~17.8]{Hardy-Wright}) that
\[
\sum_{d=1}^{\infty}\frac{\q(d)}{d^2}=\frac{\zeta(2)}{\zeta(4)}.
\]

Then,
\begin{alignat*}{2}
\sum_{k=1}^{n}\psi(k) 
&= \sum_{k=1}^{n}\Psi(k) 
= \sum_{k=1}^{n}k\sum_{d|k}\frac{\q(d)}{d} 
&&= \sum_{dd'\leq n}d'\q(d) 
= \sum_{d=1}^{n}\q(d)\sum_{d'=1}^{\floor*{\frac{n}{d}}}d' \\
&= \frac{1}{2}\sum_{d=1}^{n}\q(d)\left(\floor*{\frac{n}{d}}^2+\floor*{\frac{n}{d}}\right) 
&&= \frac{1}{2}\sum_{d=1}^{n}\q(d)\left(\frac{n^2}{d^2}+O\left(\frac{n}{d}\right)\right) \\
&= \frac{n^2}{2}\sum_{d=1}^{n}\frac{\q(d)}{d^2} + O\left(n\sum_{d=1}^{n}\frac{1}{d}\right) 
&&= \frac{n^2}{2}\sum_{d=1}^{\infty}\frac{\q(d)}{d^2} + O\left(n^2\sum_{d=n+1}^{\infty}\frac{1}{d^2}\right) + O\left(n\ln n\right) \\
&= \frac{n^2}{2}\frac{\zeta(2)}{\zeta(4)} + O\left(n\right) + O\left(n\ln n\right) 
&&= \frac{n^2}{2}\frac{\zeta(2)}{\zeta(4)} + O\left(n\ln n\right) 
\end{alignat*}
and the average order of $\psi(n)$ is $n\zeta(2)/\zeta(4)$.
This concludes the proof since, by \Cref{t:a1}, the average order of $\sigma(n)$ is $n\zeta(2)$.
\end{proof}

\begin{remark}
The result about the mean order of the Dedekind psi function can also be found in \cite{Perez}.
It is an avatar of the result in \cite[Section~18.5]{Hardy-Wright} about mean order of Euler's totient function $\varphi$, which has a similar expression to that of the Dedekind psi function $\psi$:
\[
\psi(n)=n\prod_{p|n}\left(1+\frac{1}{p}\right)
\qquad\text{and}\qquad
\varphi(n)=n\prod_{p|n}\left(1-\frac{1}{p}\right).
\]
\end{remark}

\printbibliography
\end{document}